\documentclass{l4dc2026}

\title[Closed Form HJB Solution]{Closed Form HJB Solution for Continuous-Time Optimal Control of a Non-Linear Input-Affine System}

\usepackage{times}
\usepackage{multirow}
\usepackage{cite}

\author{%
 \Name{Akash Vyas} \Email{vyas.8@iitj.ac.in}\\
 \addr Department of Mechanical Engineering, IIT Jodhpur
 \AND
 \Name{Shreyas Kumar} \Email{shreyaskumar@iisc.ac.in}\\
 \addr HiRo Lab, Cyber Physical Systems, IISc Bengaluru
 \AND
 \Name{Jayant Kumar Mohanta} \Email{jayant@iitj.ac.in}\\
 \addr Department of Mechanical Engineering, IIT Jodhpur
 \AND
 \Name{Ravi Prakash} \Email{ravipr@iisc.ac.in}\\
 \addr HiRo Lab, Cyber Physical Systems, IISc Bengaluru%
}

\begin{document}

\maketitle

\begin{abstract}%
Designing optimal controllers for nonlinear dynamical systems often relies on reinforcement learning and adaptive dynamic programming (ADP) to approximate solutions of the Hamilton Jacobi Bellman (HJB) equation. However, these methods require iterative training and depend on an initially admissible policy. This work introduces a new analytical framework that yields closed-form solutions to the HJB equation for a class of continuous-time nonlinear input-affine systems with known dynamics. Unlike ADP-based approaches, it avoids iterative learning and numerical approximation. Lyapunov theory is used to prove the asymptotic stability of the resulting closed-loop system, and theoretical guarantees are provided. The method offers a closed-form control policy derived from the HJB framework, demonstrating improved computational efficiency and optimal performance on state-of-the-art optimal control problems in the literature.%
\end{abstract}


\section{Introduction}
\label{sec:introduction}

Optimal control of dynamic systems has long attracted significant research attention, with modern approaches predominantly relying on reinforcement learning (RL) and adaptive dynamic programming (ADP) methods to handle nonlinear systems. However, this computational shift has overshadowed opportunities to obtain closed-form solutions to the Hamilton Jacobi Bellman (HJB) equation \citet{5164743,lewis2012optimal} for certain classes of input-affine systems where analytical tractability is still possible.

For linear quadratic regulator (LQR) systems with quadratic cost functions, the HJB equation reduces to the algebraic Riccati equation, which admits exact analytical solutions. In contrast, solving the HJB equation for nonlinear systems remains challenging because it requires complete knowledge of the dynamics and involves solving complex nonlinear PDEs \citet{lewis2012optimal,White1992HandbookOI}. This difficulty has motivated extensive research on approximate solution methods, particularly ADP-based approaches \citet{1239198,98e03e4c2d7743e7a9b0477e2d528af3}.

ADP methods, implemented in either model-based or model-free forms, approximate the HJB solution using value or policy iteration algorithms \citet{7362040}. Variants such as generalized policy iteration \citet{5178964} and integral policy iteration \citet{LEE2014475} aim to improve learning efficiency, while value iteration methods can converge without an initial admissible policy \citet{bian2021reinforcement,su2018online}. As a result, ADP has become a major branch of RL for both discrete-time \citet{5585774,8880690} and continuous-time systems \citet{6912014,JIANG20122699,7967695}.

Although effective, most ADP and integral reinforcement learning (IRL) methods still require partial knowledge of system dynamics (e.g., control coefficient matrices) \citet{VAMVOUDAKIS2010878,10156466}, and rely on iterative training with function approximators, which can limit practical deployment. Adaptive critic designs (ACDs) \citet{VRABIE2009477,5531586,VRABIE2009237} address this by performing online approximation of the cost function and control law, but they remain computationally demanding and lack closed-form guarantees. 



However, existing adaptive critic–based schemes typically provide only uniform ultimate boundedness (UUB) of the closed-loop trajectories. While UUB ensures that system states remain within a bounded neighborhood of the equilibrium, it does not guarantee convergence to the exact equilibrium. In safety-critical domains such as robotics or aerospace, this residual steady-state error is undesirable.

To address this limitation, this work introduces a novel drift-free augmented formulation for a class of continuous-time nonlinear input-affine systems, enabling a closed-form analytical solution of the HJB equation under known system dynamics. The proposed approach leverages Lyapunov stability theory to establish strict asymptotic convergence of the closed-loop system, thereby providing both analytical tractability and theoretical guarantees beyond those offered by adaptive critic–based methods.

\section{Mathematical Preliminaries}
\label{sec:preliminaries}
Let's consider a nonlinear input-affine dynamical system in continuous-time written as,
\begin{equation}\label{eq:control_affine_form}
\dot{x} = f(x) + g(x)\tau
\end{equation}
where, $x \in \mathbb{R}^{m}$ denotes the state, $f(x) \in \mathbb{R}^{m}$ denotes the Internal dynamics, $g(x) \in \mathbb{R}^{m \times n}$ denotes control coefficient matrix and $\tau \in \mathbb{R}^{n}$ denotes the control input. Additionally, It is assumed that the dynamics $f(x)$ and $g(x)$ are known a prior in this study.

The system is assumed to be observable and controllable, with $x=0$ as its unique equilibrium contained in a compact set $\chi$ with $f(0)=0$, ensuring that an optimal control input for \eqref{eq:control_affine_form} can be determined. To facilitate the derivation of an analytical solution, the original system \eqref{eq:control_affine_form} is rewritten in an augmented drift-free form as, 
\begin{equation}\label{eq:error_dynamics}
\dot{x}=
\begin{bmatrix}
f(x) & g(x)
\end{bmatrix}
\begin{bmatrix}
1 \\
\tau
\end{bmatrix}
=
P(x)u
\end{equation}
where, $P(x) \in \mathbb{R}^{m \times (n+1)}$, $u \in \mathbb{R}^{(n+1)}$. In the following sections, an optimal control law will be derived for this augmented system to ensure stability while minimizing a global cost function.

\subsection{The Hamilton-Jacobi-Bellman Equation}
The value function for the given system dynamics \eqref{eq:error_dynamics} over time interval $(t,\infty]$ is defined by, 
\begin{equation}
\begin{aligned}
V(x,t) &= \min_{u} J(x,u), \quad
\text{s.t.} \;  \dot{x} = f(x,u)
\end{aligned}
\label{eq:value_function}
\end{equation}

\begin{equation}
J(x,u) = \int_{t}^{\infty} \ell(x(\tau), u(\tau))\, d\tau
\label{eq:cost_function}
\end{equation}
where, $l(x,u) = \frac{1}{2}[Q(x) + u^{\top}Ru]$ is the stage cost, $Q(x_{i}) \geq 0$ and $R \in \mathbb{R}^{(n+1) \times (n+1)}$ is a positive definite matrix. A proper choice of $Q(x) = x^{\top}[Q_{0} + \gamma P(x)P(x)^{\top}]x$ where $Q_{0} \in \mathbb{R}^{m \times m}$ is symmetric, ensures that the components of the state vector are 'small' as the second term in $Q(x)$ increases the penalty on state in regions where the system is highly nonlinear and harder to control \citet{9992796}, while $R$ limits excessive control effort $u(t)$. 

The objective is to find an optimal control law $u^{*}(t)$ which minimizes the above cost function. For this we can write the Hamiltonian as, $H(x,u) = \nabla_{x}V(x,t)^{\top}\dot{x} + l(x,u)$ where, $\nabla_{x}V(x,t)^{\top}$ is the gradient of the $V(x,t)$ with respect to $x$, It is well known that the optimal control input that minimizes the cost function \eqref{eq:cost_function} also minimizes the Hamiltonian and therefore we use stationary condition to find optimal control $u^{*}(t)$, this gives the formulation popularly known as the Hamiltonian-Jacobi-Bellman (HJB) equation,
\begin{equation}\label{eq:HJB_eqn}
\min_{u} \Bigg\{\left (\frac{\partial V}{\partial x} \right)^{\top}\dot{x}(t) + l(x(t),u(t)) \Bigg\} = 0
\end{equation}
now, substituting expressions of $\dot{x}$ from \eqref{eq:error_dynamics} and $l(x(t),u(t))$ in equation \eqref{eq:HJB_eqn},
\begin{equation}\label{eq:HJB_eqn_2}
\min_{u} \Bigg\{ \left (\frac{\partial V}{\partial x} \right)^{\top}P(x)u + \frac{Q(x)}{2} + \frac{u^{\top}Ru}{2} \Bigg\} = 0
\end{equation}
Differentiating \eqref{eq:HJB_eqn_2} with respect to $u$ and setting  the derivative to 0 will give,
\begin{equation}\label{eq:optimal_u_1}
u^{*}(t) = - R^{-1}P^{\top}(x) \left (\frac{\partial V}{\partial x} \right)
\end{equation}
In order to find the expression for $\left (\frac{\partial V}{\partial x} \right)$, use the optimal control input $u^{*}(t)$ from \eqref{eq:optimal_u_1} in eqn \eqref{eq:HJB_eqn_2},
\begin{equation}\label{eq:HJB_eqn_3}
\left (\frac{\partial V}{\partial x} \right)^{\top}P(x)R^{-1}P^{\top}(x)\left (\frac{\partial V}{\partial x} \right) = Q(x)
\end{equation}
Generally to solve for the optimal control problem, equations \eqref{eq:optimal_u_1}-\eqref{eq:HJB_eqn_3} need to be solved to obtain optimal value function $\nabla_{x}V(x,t)$ first. However, finding the solution is relatively difficult due to nonlinear system \eqref{eq:control_affine_form}. Consequently, iterative methods are often employed to compute approximate optimal solutions, but they are computationally expensive, sensitive to initialization, and typically only guarantee local optimality. 

In this method we aim to construct a closed form solution for a set of nonlinear systems \eqref{eq:control_affine_form} to find an optimal control \citet{8437157}, \citet{9926505}. So, eqn \eqref{eq:HJB_eqn_3} has a quadratic form $Z^{\top}Z = z$ for $z \geq 0$ which has a solution of the form $Z = \sqrt{z}\psi$, where $\psi \in \mathbb{R}^{(n+1)}$ is a vector with $\psi^{\top}\psi=1$. Thus, we get, 
\begin{equation}\label{eq:Y}
Z = R^{-\frac{1}{2}}P^{\top}(x)\left (\frac{\partial V}{\partial x} \right) = \sqrt{Q(x)}\psi
\end{equation}
Substituting $Z$ from \eqref{eq:Y} into optimal control input $u^{*}(t)$ eqn \eqref{eq:optimal_u_1}, we get,
\begin{equation}\label{eq:optimal_u_2}
u^{*}(t) = - R^{-\frac{1}{2}}\psi\sqrt{Q(x)}
\end{equation}
The expression for $\psi$ is determined using Lyapunov Stability criteria defined as,
\begin{equation}\label{eq:lyapunov}
\mathcal{V}(x) = \frac{1}{2}x^{\top}x
\end{equation}
The time derivative of the Lyapunov function is, $\dot{\mathcal{V}}(x) = x^{\top}\dot{x} = x^{\top}P(x)u^{*}$ and substituting the control input from \eqref{eq:optimal_u_2} gives $\dot{\mathcal{V}}(x) = x^{\top}\dot{x} = -x^{\top}P(x)R^{-\frac{1}{2}}\psi\sqrt{Q(x)}$. 

Here, the choice of $\psi = \frac{P^{\top}(x)x}{\|P^{\top}(x)x\|}$ is made to directly enforce the Lyapunov stability condition while satifying the constraint $\psi^{\top}\psi=1$. This ensures that the Lyapunov derivative, $\dot{\mathcal{V}}(x)$ is negative definite.  This strategic selection is the key that enables an analytical, closed-form solution for the control law from the HJB equation. However, since $P(x)$ depends on the state, it is possible that it might lose rank or align unfavorably with $x$ at nonzero states. In such cases, the right-hand side can become arbitrarily small (or even vanish), which may locally prevent strict decrease of $\mathcal{V}(x)$. This degeneracy can give rise to saddle-type behaviors in the closed-loop dynamics.

Unlike adaptive critic–based approximations, which typically guarantee only \textit{uniform ultimate boundedness} (UUB) due to residual approximation errors, the closed-form analytical control law derived in (12) ensures a strictly negative-definite Lyapunov derivative. As a result, the proposed controller guarantees \textit{asymptotic convergence} of the system state to the equilibrium point without requiring persistent excitation or iterative parameter adaptation. \\
Hence, we get final control input as,
\begin{equation}\label{eq:optimal_u_3}
u^{*}(t) = - R^{-\frac{1}{2}}\frac{P^{\top}(x)x}{\|P^{\top}(x)x|\|}\sqrt{x^{\top}[Q_0 + \gamma P(x)P^{\top}(x)]x}
\end{equation}
where, $P(x)P^{\top}(x)=f(x)f^{\top}(x) + g(x)g^{\top}(x)$. Now, we require the state penalty $x^{\top}[Q_0 + \gamma P(x)P^{\top}(x)]x \geq 0$  $\forall x \in \chi$. So, for $x \neq 0$,
\begin{equation}\label{eq:gamma_condition_1}
 \gamma \geq - \frac{x^{\top}Q_{0}x}{\|P^{\top}(x)x\|^{2}} 
 \quad (\|P(x)^{\top} x\| \neq 0)
\end{equation}
So, parameters $\gamma$ should satisfy the condition given in \eqref{eq:gamma_condition_1} for $Q(x_{i}) \geq 0$. Now, augmented control input $u^{*}(t)$ from eqn \eqref{eq:optimal_u_3} is multiplied with $D \in \mathbb{R}^{n \times (n+1)}$ matrix,
\begin{equation}\label{eq:optimal_tau}
\tau^{*}(t)=
Du^{*}(t)=
\begin{bmatrix}
0_{n \times 1} & I_{n \times n}
\end{bmatrix}
u^{*}(t)
\end{equation}

The transformation matrix $D$ is used to extract the actual control input from the augmented optimal control $u^{*}(t)$ in \eqref{eq:optimal_u_3}. While the original system \eqref{eq:control_affine_form} uses an $n$-dimensional control input $\tau(t)$, the HJB formulation yields an $(n+1)$-dimensional solution $u^{*}(t)$, his transformation effectively removes the scalar component and retains only the $n$-dimensional control vector needed for the actual system. For practical implementation, these considerations needs to be satisfied:
\begin{enumerate}
\item Verify $\|P^{\top}(x) x\| \neq 0$ to avoid singularities in the control law.
\item Choose $\gamma$ to satisfy \eqref{eq:gamma_condition_1} over the operating region $\chi$. 
\end{enumerate}

\subsection{Proposed Optimal Tracking Control}
For trajectory tracking, the goal is to design an optimal control law $\tau^{*}(t)$ that enables the nonlinear system to follow a time-varying desired trajectory $x_{d}(t)$ with known derivative $\dot{x}_{d}(t)$ \citet{5531586}. The desired trajectory dynamics are defined as,
\begin{equation}\label{eq:desired_dynamics}
\dot{x}_d = f(x_d) + g(x)\tau_{d}
\end{equation}
where $f(x_d)$ is the internal dynamics of the system \eqref{eq:control_affine_form} expressed at the desired state $x_d$, $g(x)$ is defined in \eqref{eq:control_affine_form}, and $\tau_{d}$ is the desired control input. It is useful to note that $\dot{x}_{d}(t)$, $f(x_d)$ and $g(x)$ are known in \eqref{eq:desired_dynamics}.

Let the tracking error be defined as: $e(t) = x(t) - x_d(t)$ and its derivative with respect to time $t$ is, 
$$\dot{e}(t) = \dot{x}(t) - \dot{x}_d(t) = f(x) + g(x)\tau - f(x_d) - g(x)\tau_{d}$$
\begin{equation}\label{eq:error_dynamics_2}
\dot{e}(t) = f_{e}(e) + g(x)\tau_{e}
\end{equation}
where, $f_{e}(e) = f(x) - f(x_d)$, $\tau_{e} = \tau - \tau_{d}$. We, can convert eqn \eqref{eq:error_dynamics_2} into augmented form, 
\begin{equation}\label{eq:augmented_form}
\dot{e}=
\begin{bmatrix}
f_e(e) & g(x)
\end{bmatrix}
\begin{bmatrix}
1 \\
\tau_{e}
\end{bmatrix}
=
P_e(x)u_e
\end{equation}
To achieve optimal control of \eqref{eq:augmented_form}, the control policy $u_e$ must be chosen to minimize the infinite-horizon HJB cost function,
\begin{equation}\label{eq:cost_function_2}
V(e,t) = \min_{u_e} \left [ \int_{t}^{\infty} l(e(\tau), u_e(\tau))\, d\tau \right]
\end{equation}
where, $l(e,u_e) = \frac{1}{2}[Q(e) + u^{\top}_{e}Ru_{e}]$ is the cost function, $Q(e_{i}) \geq 0$ and $R \in \mathbb{R}^{(n+1) \times (n+1)}$ is the positive definite penalty on the control input. The Hamiltonian for the HJB tracking problem is now written as, $H(e,u_e) = \nabla_{e}V(e,t)^{\top}\dot{e} + l(e,u_e)$. Now, applying the stationary condition to it gives us the optimal control input for the tracking problem as,
\begin{equation}\label{eq:optimal_u_for_tracking}
u^{*}_{e}(t) = - R^{-1}P^{\top}_e(x) \left (\frac{\partial V}{\partial e} \right)
\end{equation}
Now, substituting \eqref{eq:optimal_u_for_tracking} into Hamiltonian yields the HJB equation for the tracking problem to be,
\begin{equation}\label{eq:hamiltonian_for_tracking}
\left (\frac{\partial V}{\partial e} \right)^{\top}P_e(x)R^{-1}P^{\top}_e(x)\left (\frac{\partial V}{\partial e} \right) = Q(e)
\end{equation}
This is a quadratic equation of the form $Z^{\top}Z = z$. It has a solution of the form $Z = \sqrt{z}\psi$, where $\psi \in \mathbb{R}^{(n+1)}$ is a vector with $\psi^{\top}\psi=1$.
and substituting the value of $Z$ in \eqref{eq:optimal_u_for_tracking}, we get,
\begin{equation}\label{eq:optimal_u_for_tracking_2}
u^{*}_e(t) = - R^{-\frac{1}{2}}\psi\sqrt{Q(e)}
\end{equation}
The expression for $\psi$ is determined using Lyapunov Stability criteria defined as,
\begin{equation}\label{eq:lyapunov}
\mathcal{V}(e) = \frac{1}{2}e^{\top}e
\end{equation}
Then the time derivative of the Lyapunov function is, $\dot{\mathcal{V}}(e) = e^{\top}\dot{e} = e^{\top}P_e(x)u^{*}_{e}$. Now, using control input expression in $\dot{\mathcal{V}}(e) = e^{\top}\dot{e} = -e^{\top}P_e(x)R^{-\frac{1}{2}}\psi\sqrt{Q(e)}$. \\
we need to choose $\psi$ such that $\dot{\mathcal{V}}(e)$ is negative definite, let $\psi = \frac{P^{\top}_{e}(x)e}{||P^{\top}_{e}(x)e||}$ it makes $\dot{\mathcal{V}}(e)$ negative definite and satisfy the constraint $\psi^{\top}\psi=1$. Hence, we get final control input as,
\begin{equation}\label{eq:optimal_u_for_tracking_3}
u^{*}_e(t) = - R^{-\frac{1}{2}}\frac{P^{\top}_e(x)e}{||P^{\top}_e(x)e||}\sqrt{e^{\top}[Q_0 + \gamma P_e(x)P_e(x)^{\top}]e}
\end{equation}
control input $\tau_{e}$ for trajectory tracking can be extracted out from \eqref{eq:optimal_u_for_tracking_3}, we can multiply it with $D \in \mathbb{R}^{n \times (n+1)}$,
\begin{equation}\label{eq:optimal_tau_for_tracking}
\tau^{*}_{e}(t)=
Du^{*}_{e}(t)=
\begin{bmatrix}
0_{n \times 1} & I_{n \times n}
\end{bmatrix}
u^{*}_e(t)
\end{equation}
Recalling $\tau_{e} = \tau - \tau_{d}$, So, eqn \eqref{eq:optimal_tau_for_tracking} can be rewritten as,
\begin{equation}\label{eq:optimal_tau_for_tracking_2}
\tau^{*}(t) = \tau^{*}_{e}(t) + \tau_{d} = Du^{*}_{e}(t) + \tau_{d}
\end{equation}
this optimal control input \eqref{eq:optimal_tau_for_tracking_2} consists of a predetermined feedforward term, $\tau_{d}$ that can be found by rearranging \eqref{eq:desired_dynamics} as,
\begin{equation}\label{eq:optimal_tau_for_tracking_3}
\tau^{*}(t) = Du^{*}_{e}(t) + g(x)^{\dagger}[\dot{x}_d - f(x_d)]
\end{equation}

\section{Results}
\label{sec:SimulationResults}
To demonstrate the effectiveness of our analytical method, we offer multiple numerical examples for set-point control and trajectory tracking with a nonlinear input-affine system on a personal computer having intel i9-13900K (13th Gen) and 64 GB of RAM was used to conduct the numerical experiments. To quantitatively assess the performance of the proposed control method, the following standard indices are considered:
\begin{enumerate}
    \item Integral of Time-weighted Squared Error (ITSE):
    It show how quickly controller drives errors to zero which highlights transient vs. steady-state trade-offs.
    \begin{equation}
        \text{ITSE} = \int_{0}^{T} te^{2}(t) \, dt
    \end{equation}

    \item Cumulative Cost: This metric is a measure of accuracy and optimal effort with smooth trajectories. 
    \begin{equation}
        J_{\tau} = \int_{0}^{T} [e^{\top}(t) e(t) + \tau^{\top}(t) \tau(t) + \dot{\tau}^{\top}(t) \dot{\tau}(t)] \, dt
    \end{equation}
    Although the methods were designed under different weighting matrices, we report a normalized evaluation cost using $Q=I$, $R=I$ for fair comparison.  

    \item Computation Time:
    The wall-clock time required to reduce the tracking error below $10^{-3}$ is reported; if this threshold is not reached, the time for the full horizon $[0,T]$ is used.
\end{enumerate}

\subsection{Example I}
Consider this continuous-time nonlinear system from \citet{5531586} and \citet{VAMVOUDAKIS2010878},
\begin{multline}\label{eq:example_1}
\begin{bmatrix}
\dot{x}_1 \\
\dot{x}_2
\end{bmatrix}
=
\begin{bmatrix}
-x_1 + x_2 \\
-\frac{x_1}{2} - \frac{x_2\big(1-(\cos(2x_1)+2)^2\big)}{2} 
\end{bmatrix}
+
\begin{bmatrix}
0 \\
\cos(2x_1)+2
\end{bmatrix}
\tau
\end{multline}
\begin{figure*}[!htb]
    \centering
    \includegraphics[width=\textwidth]{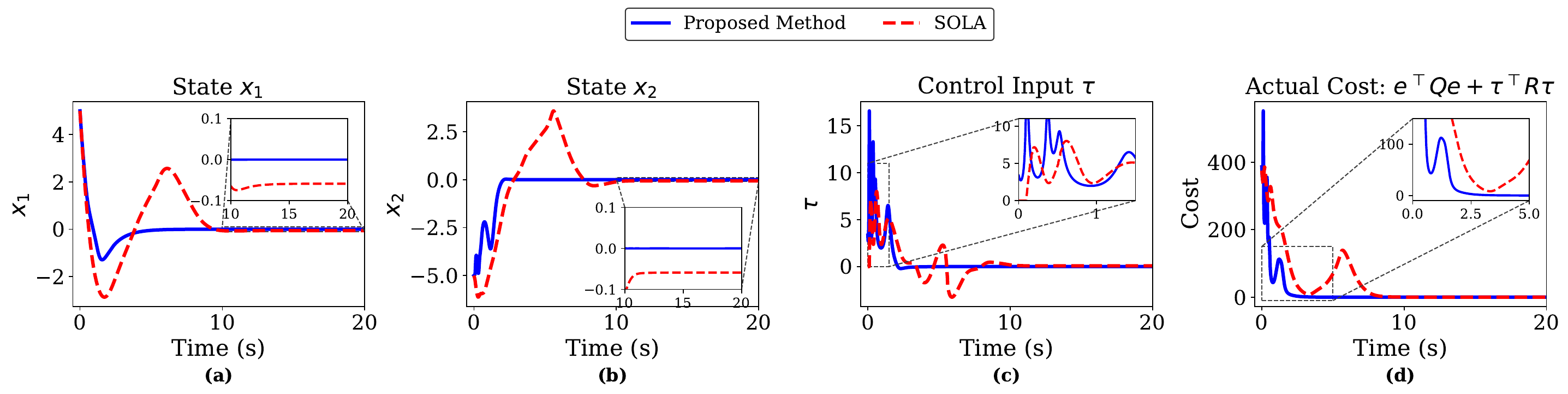}
    \caption{System response: (a) State variables $x_1$, (b) State variables $x_2$, (c) Control input $\tau$, (d) Cost function}
    \label{fig:figure_1}
\end{figure*}

We conducted a quantitative comparison study of the proposed method with a single online approximator (SOLA) based scheme \citet{5531586}. For SOLA, the stage cost was $l(x,u) = Q(x) + u^{\top}Ru$ with $Q(x) = x^{\top}x$, $R=1$ and the basis vector 
\[
\phi(x) = [x_{1}, \,x_{2}, \,x_{1}x_{2}, \,x_{1}^2, \,x_{2}^2, \,x_{1}^2cos(2x_{1})^2, \,x_{1}^3]^{\top}
\]
was used. The tuning parameters $\alpha_1=25.0$, $\alpha_2=0.01$, and the initialization of all NN-weights to zero were taken directly from \citet{5531586}. Initial condition for both methods was $x(0)=[5, -5]^{\top}$ and for the proposed method $Q_0 = I_{2 \times 2}$, $R=I_{2 \times 2}$ and $\gamma = 1.0$ were used.\\

Fig.~\ref{fig:figure_1} illustrates the comparative performance over $t \in [0, T]$. Subplots (a-b) shows the state convergence to the equilibrium point $x=0$ for both methods while subplots (c-d) illustrates the control effort and cost metric. For a fair comparison, the controllers were simulated without adding any probing noise. While persistent excitation is generally required for parameter convergence in adaptive control, practical control objectives can often be achieved without artificial excitation, as demonstrated by successful convergence in the results. 
\begin{table}[htbp]
\caption{Performance Comparison of Optimal Control Methods for Example-I}
\centering
\setlength{\tabcolsep}{2pt} 
\begin{tabular}{|c|c|c|c|c|}
\hline
\textbf{Method} & \textbf{ITSE} & \textbf{Cumulative Cost} & \textbf{Computation Time (s)} \\  
\hline
HJB-SOLA \citet{5531586} & 4615.71 & 1546.702 & 2.235 \\
\hline
Proposed method & 35.977 & 876.785 & 0.201 \\
\hline
\end{tabular}
\label{tab:comparison_1}
\end{table}
Table~\ref{tab:comparison_1} shows that the proposed HJB-based method achieves the lowest ITSE and cumulative cost, indicating faster convergence with lower control effort than the SOLA baseline~\citet{5531586}. It also records the shortest computation time, demonstrating computational efficiency.

\subsection{Example II}
\begin{figure*}[!htb]
    \centering
    \includegraphics[width=\textwidth]{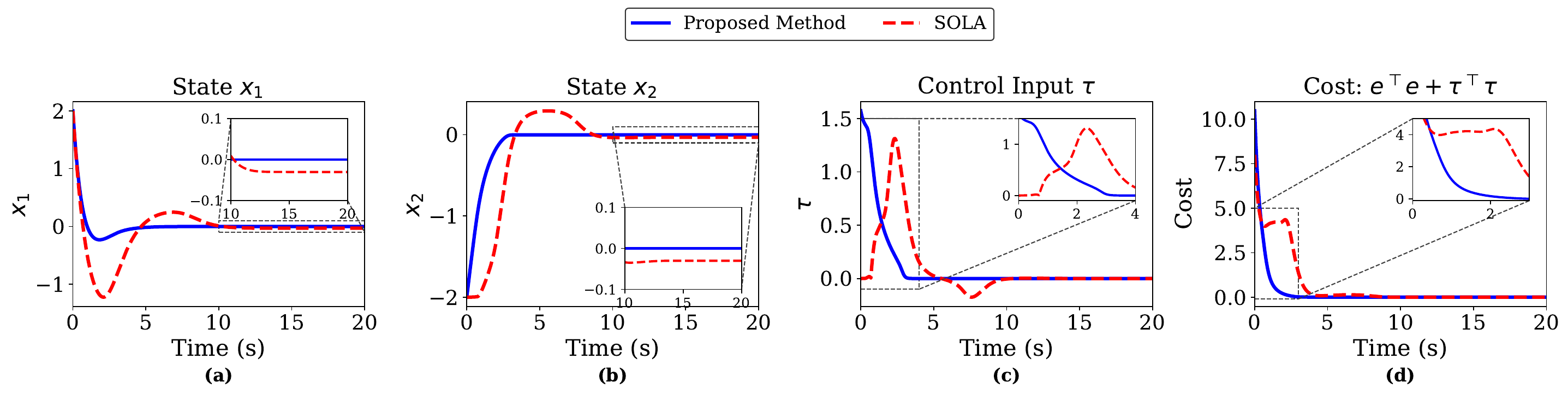}
    \includegraphics[width=\textwidth]{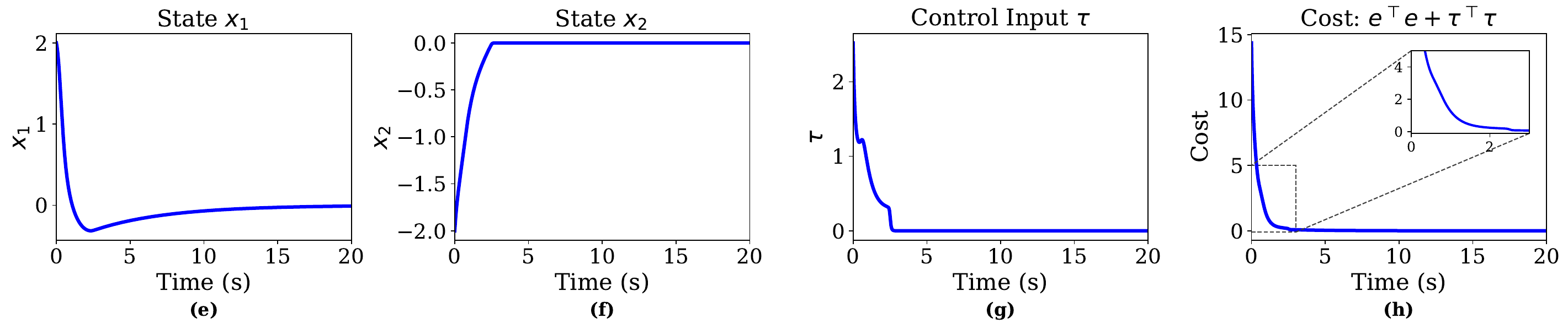}
    \caption{System response: (a) State variables $x_1$ for Case-1, (b) State variables $x_2$ for Case-1, (c) Control input $\tau$ for Case-1, (d) cost function for Case-1, (e) State variables $x_1$ for Case-2, (f) State variables $x_2$ for Case-2, (g) Control input $\tau$ for Case-2, (h) cost function for Case-2}
    \label{fig:figure_2}
\end{figure*}

The following nonlinear system studied here is taken from \citet{Lin01012000} and \citet{6707098},
\begin{equation}\label{eq:example_2}
\begin{bmatrix}
\dot{x}_1 \\
\dot{x}_2
\end{bmatrix}
=
\begin{bmatrix}
x_2 + \lambda_1x_1cos\left(\frac{1}{x_2+\lambda_2} \right) + \lambda_3x_2sin(\lambda_4x_1x_2) \\
0
\end{bmatrix}
+
\begin{bmatrix}
0 \\
1
\end{bmatrix}
\tau
\end{equation}
Now, we simulate the close-loop system based on optimal control policy \eqref{eq:optimal_tau}. Set the initial state be $x(0)=[2,-2]^{\top}$ and the simulation results for the following three cases are given in Fig.~\ref{fig:figure_2}, for the proposed method where $Q_0 = I_{2 \times 2}$, $R=I_{2 \times 2}$, $\gamma_{case-1} = 0.5$, $\gamma_{case-2} = 0.5$ and $\gamma_{case-3} = 1.0$ respectively.
\begin{enumerate}
\item Case 1: $\lambda_1 = -1$, $\lambda_2 = -100$, $\lambda_3 = 0$, $\lambda_4 = -100$
\item Case 2: $\lambda_1 = -0.2$, $\lambda_2 = 100$, $\lambda_3 = 1$, $\lambda_4 = -1$
\end{enumerate}

From the results shown in Fig.~\ref{fig:figure_2}, it can be observed that under the optimal control policy \eqref{eq:optimal_tau}, system states converge to their respective equilibrium points: subplots (a–d) correspond to Case-1, subplots (e–h) to Case-2, and subplots (i–l) to Case-3. For a fair evaluation, all methods were simulated over the time interval $t \in [0, T]$ without adding any probing noise, and the basis functions were kept the same as those used in \citet{5531586}. Moreover, no explicit tuning of $Q$ and $R$ was performed for any method; their values were fixed at 1. 
\begin{table}[htbp]
\caption{Performance Comparison of Optimal Control Methods for Example-II}
\centering
\setlength{\tabcolsep}{2pt} 
\begin{tabular}{|c|c|c|c|c|}
\hline
\textbf{Case} & \textbf{Method} & \textbf{ITSE} & \textbf{Cumulative Cost} & \textbf{Computation Time (s)}\\
\hline

\multirow{2}{*}{Case 1}
 & HJB-SOLA \citet{5531586} & 40.59 & 16.62 & 2.085\\ \cline{2-5}
 & Proposed method & 2.036 & 6.097 & 0.191\\ 
 \hline

\multirow{2}{*}{Case 2}
 & HJB-SOLA \citet{5531586} & N/C & N/C & N/C\\ \cline{2-5}
 & Proposed method & 2.684 & 14.859 & 0.649\\ 
 \hline

\end{tabular}\\
\vspace{3pt}
\scriptsize{$^{*}$N/C: Not Converged}
\label{tab:comparison_2}
\end{table}
The results of the comparative study is presented in Table~\ref{tab:comparison_2}, the proposed HJB-based method achieves lowest ITSE, demonstrating superior transient response and faster error decay compared to \citet{5531586}. 
\begin{figure*}[!htb]
    \centering
    \includegraphics[width=\textwidth]{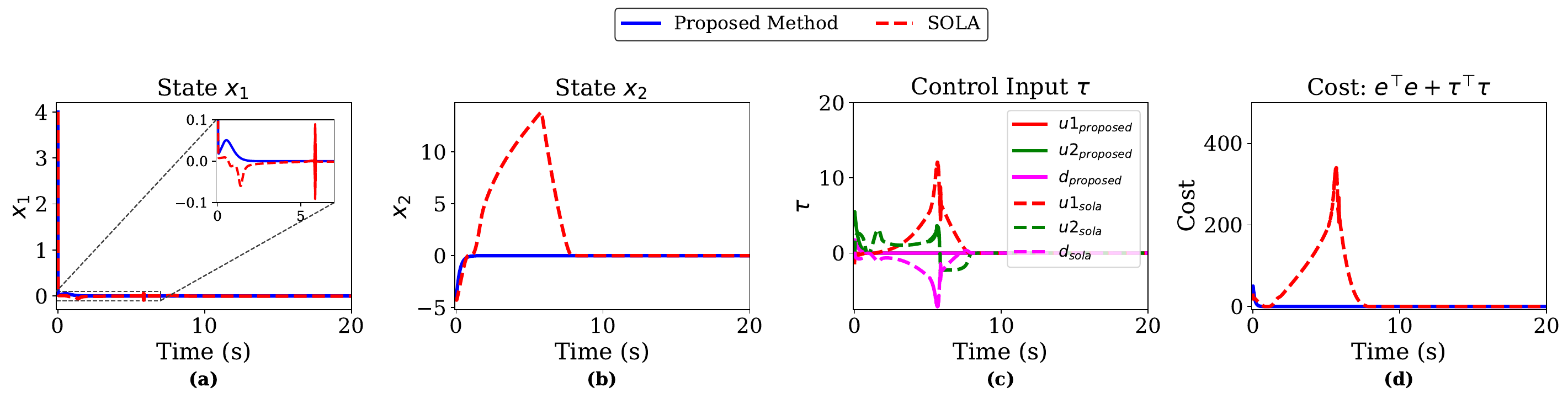}
    \caption{System response for Case 3: (a) State variables $x_1$, (b) State variables $x_2$, (c) Control input $\tau$, (d) Cost function}
    \label{fig:figure_5}
\end{figure*}
\subsection{Example III}
In this example, we are trying to solve a different problem from differential game theory and this following nonlinear system studied here is taken from \citet{5717676}. Let's consider the nonlinear system in the presence of disturbances as,
\begin{equation}
\dot{x} = f(x) + g(x)u + k(x)d 
\end{equation}
where, $x \in \mathbb{R}^{m}$ denotes the state, $u \in \mathbb{R}^{n}$ denotes the control input, $d \in \mathbb{R}^{w}$ denotes the disturbance, $k(x) \in \mathbb{R}^{m \times w}$ is bounded according to $\|k(x)\|_{F} \leq k_{M}$ for a known constant $k_M$, where $\|.\|_{F}$ is the Frobenius norm. These systems are two-player zero-sum differential game, governed by the Hamilton Jacobi Isaacs (HJI) equation,
\begin{multline}\label{eq:example_3}
\begin{bmatrix}
\dot{x}_1 \\
\dot{x}_2
\end{bmatrix}
=
\begin{bmatrix}
-\frac{29x_{1} + 87x_{1}x^{2}_{2}}{8} - \frac{2x_{2} + 3x_{2}x^{2}_{1}}{4} \\
-\frac{x_{1} + 3x_{1}x^{2}_{2}}{4}
\end{bmatrix}
+
\begin{bmatrix}
1 & 0\\
0 & 3
\end{bmatrix}
\begin{bmatrix}
u_1 \\
u_2
\end{bmatrix}
\\[4pt]
+
\begin{bmatrix}
0.5\\
1
\end{bmatrix} d
\end{multline}

We simulated the close-loop system under the optimal control policy \eqref{eq:optimal_tau}. For this, the control input was assumed to be a vector $\tau = [u_1, u_2, d]^{\top}$ and band the initial state was set as $x(0)=[4,-4]^{\top}$ with $Q_0 = I_{2 \times 2}$, $R=I_{3 \times 3}$ and $\gamma = 0.1$ respectively for the proposed method. \\
From the results shown in Fig.~\ref{fig:figure_5}(a)-(b), it can be observed that under the action of the optimal control policy \eqref{eq:optimal_tau}, the states of system \eqref{eq:example_3} can converge to the equilibrium point.

We also conducted a quantitative comparison of the proposed method with a single online approximator (SOLA) based scheme reported in \citet{5717676}. For the HJI problem, the penalties were set as $Q(x) = 2((2x_{1} + 6x_{1}x^{2}_{2})^{2} + (4x_{2} + 6x^{2}_{1}x_{2})^{2})$ and $R=1$. The basis vector $\phi(x) = [x_{1}, x_{2}, x_{1}x_{2}, x_{1}^2, x_{2}^2, x_{1}^{2}x_{2}^{2}, x_{1}^3]^{\top}$ was used. The tuning parameters were chosen as $\alpha_1=200.0$, $\alpha_2=0.01$. The initial condition was set to $x(0)=[4, -4]^{\top}$ and all neural network weights were initialized to zero following \citet{5717676}. \\
\begin{table}[htbp]
\caption{Performance Comparison of Optimal Control Methods for Example-III}
\centering
\setlength{\tabcolsep}{2pt} 
\begin{tabular}{|c|c|c|c|c|}
\hline
\textbf{Method} & \textbf{ITSE} & \textbf{Cumulative Cost} & \textbf{Computation Time (s)}\\
\hline
HJI-SOLA \citet{5717676} & 9567621.127 & 25233371.564 & 0.561\\
\hline
Proposed method & 1.155 & 979.797 & 0.0705\\
\hline
\end{tabular}
\label{tab:comparison_3}
\end{table}
For a fair evaluation, both methods were simulated for time $t \in [0, T]$ without adding any probing noise and system \eqref{eq:example_3} was able to converge without it. Table~\ref{tab:comparison_3} shows that the proposed HJB-based method achieves the lowest ITSE and cumulative cost, indicating faster convergence with lower control effort than \citet{5717676}. It also records the shortest computation time, demonstrating its computational efficiency.

\section{Conclusion}
 In this paper, a novel analytical solution was presented for the Hamilton Jacobi Bellman equation applied to continuous-time nonlinear affine systems. Unlike popular ADP-based methods, the proposed approach is purely analytical, does not require training, and is computationally efficient. For systems with known dynamics, both optimal regulation and tracking control problems were addressed, and Lyapunov techniques were employed to demonstrate stability of the resulting optimal control scheme. Extensive numerical tests on multiple benchmark problems illustrate the advantages of the proposed method. Proposed optimal tracking controller also exhibited satisfactory performance, though results are omitted due to space constraints.

\bibliography{ref}

\end{document}